\documentclass[11pt]{article}

\usepackage{amssymb}
\usepackage{amsthm}
\usepackage[english]{babel}
\usepackage{mathtools}
\usepackage[all,cmtip]{xy}
\usepackage{tikz}
\usetikzlibrary{matrix,arrows}

\theoremstyle{plain}
\newtheorem{theorem}{Theorem}
\newtheorem{lemma}[theorem]{Lemma}

\theoremstyle{definition}

\theoremstyle{remark}

\title{ Quasi-projective posets, lattices, permutations, graphs, digraphs, hypergraphs, point-line geometries }

\author{\' Eva Jung\'abel}

\date{}

\begin{document}

\maketitle

\begin{abstract} A structure $\cal S$ is quasi-projective if for every structure $\cal T$, for every homomorphism $f : {\cal S} \rightarrow {\cal T}$ and every epimorphism $j: {\cal S}\rightarrow {\cal T}$ there is an endomorphism $\phi$ of $\cal S$ such that $\phi\circ j=f$. In this paper, we characterise the quasi-projective posets and lattices of arbitrary cardinalities, finite permutations, graphs and digraphs of arbitrary cardinalities with loops and without loops, finite hypergraphs, and finite point-line geometries.

\end{abstract}

\let\thefootnote\relax

\section{Introduction}\label{intro}

By a structure we mean a set together with an indexed set of relations and operations on it. A first-order structure $\cal A$ is called {\it homogeneous} if any isomorphism between two finitely generated substructures of $\cal A$ is induced by some automorphism of $\cal A$. In several classes of combinatorial structures the homogeneous structures are classified. 

P. Cameron and J. Ne\v set\v ril \cite{cameron-nesetril} introduced the following variant of homogeneity: a structure is called {\it homomorphism-homogeneous} if every homomorphism between finite induced substructures extends to an endomorphism of the structure. 

Ho\-mo\-mor\-phism-ho\-mo\-ge\-neous graphs were investigated by Rusinov and Schweitzer  in \cite{rusinov-schweitzer}. Among others it is shown that the problem of deciding if the graph is homomorhism-homogeneous is coNP-complete. Finite algebras also harbour some classes of high computational complexity \cite{masulovic3}, hence we cannot expect a brief classification in case of algebraic structures in general. A characterisation of all ho\-mo\-mor\-phism-ho\-mo\-ge\-neous partial orders of arbitrary cardinalities with non-strict relation is given  by Ma\v sulovi\' c  \cite{masulovic2} and {Cameron and Lockett} \cite{cameron-lockett}, independently. Several  other ho\-mo\-mor\-phism-ho\-mo\-ge\-neous structures are characterised including finite ho\-mo\-mor\-phism-homogeneous permutations by Dolinka and Jung\'abel \cite{dolinka-jungabel}. In \cite{jungabel1} it is shown that a the category of point-line geometries is equivalent to a certain subclass of 3-uniform hypergraphs and therefore a first step in cha\-rac\-teri\-sing homomorphism-homogeneous hypergraphs could be the investigation of point-line geometries. In \cite{jungabel1} authors show the local be\-ha\-vior of finite homomorphism-homogeneous point-line geometries with $k$ non-intersecting regular lines where every point lies on a regular line and they give a des\-crip\-tion of  these ho\-mo\-mor\-phism-ho\-mo\-ge\-neous finite point-line geometries with two non-inter\-sec\-ting regular lines. A characterisation of finite ho\-mo\-mor\-phism-homogeneous point-line geometries containing two regular intersecting lines  is given by Ma\v sulovi\' c \cite{masulovic1}.

\vskip3mm
Interestingly, the concept of homomorphism-homogeneity is not a recent concept, it exists under the name of {\it quasi-injectivity} just with the slight difference. A structure is said to be quasi-injective if every homomorphism from an arbitrary substructure of the structure into the structure extends to an endomorphism of the structure.

All quasi-injective Abelian groups are described \cite{fuchs1} as finite quasi-injective groups \cite{bertholf-walls}. Infinite quasi-injective groups \cite{tomkinson} are partly characterised. There are results about quasi-injective modules by Johnson and Wong \cite{johnson-wong1}, Harada \cite{harada1}, Faith and Utumi \cite{faith-utumi1}, Fuchs \cite{fuchs2} and others.

The dual concept of quasi-injectivity is {\it quasi-projectivity} and it was introduced for modules by Wu and Jans \cite{wu-jans1} in 1967. Phrased in terms of diagrams, the module $\cal S$ is quasi-injective if every diagram

\centerline{
\xymatrix@=5em{
  0 \ar[r] 
& {\cal T} \ar[r]^j \ar[d]^f
&{\cal S}\\
& {\cal S}}
}

\noindent
can be embedded in a commutative diagram

\centerline{
\xymatrix@=5em{
  0 \ar[r] 
& {\cal T} \ar[r]^j \ar[d]^f
&{\cal S} \ar@{-->}[dl]_{\phi} \\
& {\cal S}}
}

\noindent
where $\cal T$ is a submodule of $\cal S$, $j$ is the monomorphism  and $f$ is a homomorphism of $\cal T$ into $\cal S$.

The module $\cal S$ is said to be quasi-projective if every diagram

\centerline{
\xymatrix@=5em{
& {\cal S} \ar[d]^f \\
{\cal S} \ar[r]^j
& {\cal S/T} \ar[r]
& 0}
}

\noindent
can be embedded in a commutative diagram

\centerline{
\xymatrix@=5em{
& {\cal S} \ar[d]^f \ar@{-->}[dl]_-{\phi} \\
{\cal S} \ar[r]^j
& {\cal S/T} \ar[r]
& 0}
}

\noindent
where $\cal T$ is a submodule of $\cal S$, $j$ is the epimorphism  and $f$ is a homomorphism of $\cal T$ into $\cal S$.

In \cite{wu-jans1} some properties of quasi-projective modules are shown, a structure theorem for indecomposable finitely generated quasi-projectives over semi-perfect rings is obtained and the finitely generated quasi-projective Abelian groups are described. The general case of quasi-projective Abelian groups is characterised by Fuchs and Rangaswamy \cite{fuchs-rangaswamy1}. A decomposition theorem that is a characterisation for quasi-projective modules over left perfect rings is given by Koehler \cite{koehler}.

We say that a structure $\cal S$ is quasi-projective if for every structure $\cal T$, for every homomorphism $f : {\cal S} \rightarrow {\cal T}$ and every epimorphism $j: {\cal S}\rightarrow {\cal T}$ there is an endomorphism $\phi$ of $\cal S$ such that $\phi\circ j=f$. In this paper we characterise quasi-projective posets and lattices of arbitrary cardinalities, finite permutations, graphs and digraphs of arbitrary cardinalities  with loops and without loops, finite hypergraphs and finite point-line geometries.

\section{Quasi-projective partially ordered sets and lattices}\label{qpposets}

The purpose of this section is to characterise quasi-projective partially ordered sets and lattices of arbitrary cardinalities.

A (non-strict) partial order is a binary relation $\leqslant$ over a set $L$ satisfying the following axioms

\begin{enumerate}
\item[$\bullet$]reflexivity: $a\leqslant a$ for all $a\in L$ ,
\item[$\bullet$] antisymmetry: if $a\leqslant b$ and $b\leqslant a$, then $a=b$,
\item[$\bullet$] transitivity: if $a\leqslant b$ and $b\leqslant c$, then $a\leqslant c$.
\end{enumerate}

\noindent
If for any $a$ and $b$, either $a\leqslant b$ or $b\leqslant a$, then it is called a total linear order.

\vskip3mm
A relation $<$ is a strict order on $L$. We write $a<b$ if $a\leqslant b$ and $a\not= b$.

Let $(L,\leqslant)$ be a partial ordered set. For $\emptyset \neq K \subseteq L$, a {\it substructure} of $(L,\leqslant)$ generated by $K$ is a partially ordered set $(K, \leqslant_{K})$ such that $\leqslant_K$ is a $\leqslant \cap K^2$ relation on $K$. 

$(L,\leqslant)$ is called {\it connected} partial ordered set if for every pair $a,b\in L$, there is a finite sequence $a=c_1,c_2,…,c_n=b$, with each $c_i\in L$, such that $c_i\leqslant c_{i+1}$ or $c_i\geqslant c_{i+1}$ for each $i=1,2,\ldots,n-1$. A {\it connected component} in a partial ordered set $(L,\leqslant)$ is a maximal connected substructure.

For a partial ordered set $(L,\leqslant)$, a {\it chain}, $C$, is a structure such that every pair of elements of $C$ is comparable. If no elements are comparable, then we say that the partially ordered set is an antichain.

A {\it lattice} $(L,\leqslant)$ is a partially ordered set in which every two elements have a unique a least upper bound and a unique greatest lower bound.

A homomorphism $f : {\cal A} \rightarrow {\cal B}$ of relational structures of the same similarity type $\tau$ is a mapping between their underlying sets preserving their relations in the sense that for each $n$-ary relational symbol $R \in \tau$ we have that $(a_1, \ldots , a_n) \in R^{\cal A}$ implies $(f(a_1), \ldots , f(a_n)) \in R^{\cal B}$. We add the adjective ‘partial’ to these morphisms whenever the domain of $f$ is restricted to some (induced) substructure of $\cal A$.

\begin{lemma}\label{p1} 
A quasi-projective partially ordered set is either connected or an antichain.
\end{lemma}

\begin{proof} Let $K_1$ and $K_2$ be two different connected components in a quasi-projective partially ordered set $(L,\leqslant)$ where there are three different elements $u$, $v$ and $w$ such that $u, v\in K_1$, $u< v$, and $w \in K_2$. Such elements exist, because $L$ is not connected and it is not an antichain. Let $T=\{a,b\}$, where $a \leqslant b$. We define $f$ and $j$ in the following way:

$$f(x)=\begin{cases}    b & if \hbox{ } x>u \\ a & otherwise \end{cases}$$

$$j(x)=\begin{cases}  a & if \hbox{ } x\in K_1 \\  b & otherwise \end{cases}.$$

We claim that $f$ is a homomorphism, because if $x\leqslant u$ or $x>u$, then we have $a\leqslant a$ or $b>a$. If $x$ and $u$ are incomparable elements, but they are in the same connected component, then all elements which are incomparable with $u$ are mapped to $a$, other elements to $a$ or $b$. The situation where there is an incomparable element $y$ with $u$ and an element $z$ such that $z>u$ and $z<y$ cannot happen, because then $u<y$. One can easily see that $j$ is an epimorphism which can be seen very directly.

Assume $\phi$ is a homomorphism with $j \circ \phi=f$.  We have $f(u)=a$ and $f(v)=b$, the preimage of $a$ is $K_1$ by the map $j$ and $j^{-1}(b)\notin K_1$. So, if there is $\phi$ such that $f=j \circ \phi$, then $\phi(u)\in K_1$ and $\phi(v)\notin K_1$ contradicting the fact that $u$ and $v$ are comparable.

\end{proof}

\begin{theorem}\label{p2}
A partially ordered set is quasi-projective if and only if it is a chain or an antichain.
\end{theorem}

\begin{proof}

Let $(L,\leqslant)$ be a partially ordered set. Suppose that it is not a chain or an antichain. From Lemma \ref{p1} we know that it is connected. Let $u, v$ and $w$ be three different elements such that $u< v$, $u< w$ and $v$ and $w$ are incomparable. The situation is analogous when there different are elements $u, v$ and $w$ such that $v< u$, $w< u$ and $v$ and $w$ are incomparable. 

Suppose there is an element $z$ such that $v< z$, then let $T=\{a,b,c,d\}$ be a chain. Otherwise, if there does not exist an element $z$ such that $v< z$, then let $T=\{a,b,c\}$ be an another chain. We define $f$ and $j$ in the following way:

$$f(x)=\begin{cases}   c & if \hbox{ } x>u \\ b & otherwise \end{cases}$$

$$j(x)=\begin{cases}  b & if \hbox{ } x=v \\ a & if \hbox{ } x< v \\ d & if \hbox{ } x>v \\ c & otherwise \end{cases}$$

One can verify in the same manner as in the proof of Lemma \ref{p1} that $f$ is a homomorphism. The function $j$ is epimorphism, because there is no $x$ such that $x<v$ and $w<x$, because elements $v$ and $w$ are incomparable. Also, there is no $x$ such that $v<x$ and $x< w$. Also, if $y$ and $v$ are incomparable elements, then all elements which are incomparable with $v$ are mapped to $c$, other elements to $a$, $b$ or $d$. The situation where there is an incomparable element $y$ with $v$ and an element $z$ such that $z>v$ and $z<y$ cannot happen, because then $v<y$. Also, there is not an element $z$ such that $z<v$ and $z>y$, because then $v>y$.

Now we prove, there is no homomorphism $\phi$ such that $j \circ \phi=f$.  We have $f(u)=b$ and $f(v)=c$, the preimage of $b$ is $\{v\}$ by the map $j$ and the preimage of $c$ by the map j are incomparable elements with $v$. So, if there is $\phi$ such that $f=j \circ \phi$, then $\phi(u)=v$ and every element from $\phi(v)$ is an incomparable element with $v$. This contradicts the fact that $u$ and $v$ are comparable.

Let $(L,\leqslant)$ be a partially ordered set where it is a chain or an antichain. Any function from an antichain is a homomorphism, so these are quasi-projective partially ordered sets. 

The image of a surjective homomorphism  from a chain is also a chain. We choose an element form the set $j^{-1}(f(l))$, where $l\in L$, and put $\phi(l)=j^{-1}(f(l))$. We fix the element $\phi(l)$. For every other element $l_1\not=l$, $l_1\in L$, such that $f(l_1)=f(l)$, we choose $\phi(l)$ and put $\phi(l_1)=\phi(l)$, Thus it is easy to verify that if $\phi(x) \in j^{-1}(f(x))$ for every $x \in L$, then $\phi$ is a homomorphism with $j \circ \phi=f$.

\end{proof}

\begin{theorem}\label{p3}
 The lattice $(L,\leqslant)$ is quasi-projective if and only if it is a chain.
\end{theorem}

\section{Quasi-projective permutations}\label{qpperm}

Here we characterise quasi-projective finite permutations.

Recall \cite{cameron1} that a permutation is a relational structure $\pi = (A,<_1,<_2)$, where the underlying set $A$ is equipped with two (strict) linear orders $<_1$ and $<_2$. For example, the permutation on $A = \{e, h, i, t, w\}$ represented by the sequence black would correspond to a pair of chains on $A$ where $e <_1 h <_1 i <_1 t <_1 w$ (the basic alphabetic order) and $w <_2 h <_2 i <_2 t <_2 e$.

We define conceptually no different a permutation on $A$ to be a structure $\pi = (A,\leqslant_1,\leqslant_2)$, where $\leqslant_1$ and $\leqslant_2$ are two reflexive total orders of $A$.

We define homomorphism in the same way as in Section \ref{posets}.

\begin{theorem}\label{pr1}
Every ${\pi}$ permutation is quasi-projective.
\end{theorem}

\begin{proof}

Let ${\pi}=(A,\leqslant_1,\leqslant_2)$ be a permutation. Let $f$ be a homomorphism from ${\pi}$ to a permutation ${\sigma}$ and $j$ be an epimorphism from ${\pi}$ to ${\sigma}$. Let $a\in A$ be arbitrary. Because the function $j$ is an epimorphism, $f(a)\in j(A)$. If there are more than one $f(a)$ in $j(A)$, then let it be $j(b)=f(a)$ for an element $b\in A$. Because $j$ is an epimorphism, there does not exist an element $c$ such that $a\leqslant_1 c \leqslant_1 b$ and $j(c)\not=f(a)$. So, preimages of $f(a)$, $a\in A$ , by the map $j$ are chains with respect to the relation $\leqslant_1$. This is also true for the realtion $\leqslant_2$.

For every element $a\in A$ we choose the least element in the chain $j^{-1}(f(a))$ with respect to the relation $\leqslant_1$ and put: $\phi(a)=\min_{\leqslant_1} j^{-1}(f(a))$. It is homomorphism, because suppose that $a\leqslant_1 b$ and $a\leqslant_2 b$. Then $f(a)\leq_1 f(b)$ and $f(a)\leq_2 f(b)$. Also, we have $\min_{\leqslant_1} j^{-1}(f(a)) \leqslant_1 \min_{\leqslant_1} j^{-1}(f(b))$ and $\min_{\leqslant_1} j^{-1}(f(a)) \leq_2 \min_{\leqslant_1} j^{-1}(f(b))$, otherwise $j$ is not an epimorphism. So, we have $\phi(a)\leqslant_1 \phi(b)$ and $\phi(a)\leqslant_2 \phi(b)$. It is analogous when $a\leqslant_1 b$ and $b\leqslant_2 a$.

\end{proof}

\section{Quasi-projective graphs and digraphs with or without loops}\label{qpgad}

In this section first we characterise quasi-projective graphs and digraphs of arbitrary cardinalities when loops are not allowed, then when loops are allowed.

A graph is a pair ${\cal G} = (V, E)$, where $V$ is a set whose elements are called {\it vertices} and $E$ is a set of paired vertices, whose elements are called {\it edges}. A loop is an edge that join a vertex to itself.

For $\emptyset \neq H \subseteq G$, a {\it subgraph} of ${\cal G}=(V, E)$ generated by $H$ is a graph $(H, E\cap H^2)$. 

If $E=\emptyset$ we say that a graph ${\cal G}=(V,E)$ is empty graph. $K_n=(V,E)$ is called a complete graph if $|V|=n$ and $E=V^2$ if loops are allowed, unless $E=V^2\setminus\{(x,x)\mid x\in V\}$. We say that a graph is e-empty if there are not edges between vertices (loops can be).

A graph $(V,E)$ is said to be {\it connected} if for every pair $a,b\in V$, there is a finite sequence $a=c_1,c_2,…,c_n=b$, with each $c_i\in V$, such that $(c_i,c_{i+1})\in E$ or $(c_{i+1}, c_i)\in E$ for each $i=1,2,…,n-1$. A {\it connected component} in a graph $(V,E)$ is a maximal connected substructure.

A {\it directed graph} is an ordered pair ${\cal D} = (V, A)$ where $V$ is a set whose elements are called {\it vertices}; $A\subseteq \{(x,y)\mid (x,y)\in V^{2}\}$ is a set of ordered pairs of vertices, called {\it arrows}, {\it directed edges}. If loops are not allowed, then $A\subseteq \{(x,y)\mid (x,y)\in V^{2}\;{\textrm {and}}\;x\neq y\}$.

For $\emptyset \neq H \subseteq D$, a {\it subdigraph} of ${\cal D}=(V, A)$ generated by $H$ is a digraph $(H, A\cap H^2)$. 

$D_n=(V,A)$ is called a complete digraph if $|V|=n$ and $A=V^2$ if loops are allowed, unless $A=V^2\setminus\{(x,x)\mid x\in V\}$. We say that a digraph is e-empty if there are not arrows between vertices (loops can be).

We define homomorphisms in the same way as in Section \ref{qpposets}.

\begin{lemma}\label{g1}
A quasi-projective graph, where loops are not allowed, is a complete or an empty graph.
\end{lemma}

\begin{proof} Let a graph ${\cal G}=(V,E)$ be quasi-projective where loops are not allowed and let $V=\{v_1, v_2,\ldots, v_n\}$. Suppose to the contrary that it is neither an empty nor a complete graph. Let  $v_i,v_j\in V$ be two vertices such that $(v_i,v_j)\notin E$. Let $(v_k,v_l)\in E$, $v_k, v_l\in V$. We can permutate the indices such that $(v_1,v_2)\notin E$ and $(v_3,v_4)\in E$ or $(v_1,v_2)\notin E$ and $(v_2,v_3)\in E$. Without loss of generality, we suppose the previous case, while the proof for the latter is analogous. Let $T=K_n=\{k_1,k_2,\ldots,k_n\}$. We define $f$ and $j$ in the following way:

$$f(x)=\begin{cases}   k_i & if \hbox{ } x=v_i \end{cases},$$

$$j(x)=\begin{cases}  k_{i+2\pmod{n}} & if \hbox{ } x=v_i \end{cases}.$$

The map $f$ is homomorphism and the map $j$ is epimorphism obviously. Now we prove, there is no homomorphism $\phi$ such that $j \circ \phi=f$.  We have $f(v_3)=k_3$ and $f(v_4)=k_4$, the preimage of $k_3$ is $\{v_1\}$ by the map $j$ and $j^{-1}(k_4)=\{v_2\}$. So, if there is $\phi$ such that $f=j \circ \phi$, then $\phi(v_3)=v_1$ and $\phi(v_4)=v_2$. This contradicts the fact that $v_1, v_2\notin E$.
\end{proof}

\begin{lemma}\label{g2}
A complete or an empty graph, where loops are not allowed, is quasi-projective.
\end{lemma}

\begin{proof} From Lemma \ref{g1} we know that if a graph, where loops are not allowed, is quasi-projective, then it is a complete or an empty graph.

Now, let ${\cal G}=(V,E)$ be a complete or an empty graph where loops are not allowed. Every function of an empty graph is a homomorphism, so these are quasi-projective graphs. A homomorphism of a complete graph (to any loopless graph) is injective. Thus $j$ is an isomorphism. Thus $\phi=j^{-1}f$ satisfies the required condition.
\end{proof}

\begin{theorem}\label{g3}
A digraph, where loops are not allowed, is quasi-projective if and only it is a complete or an empty graph.
\end{theorem}

\begin{proof} Let  ${\cal D}=(V,A)$ be a complete or an empty digraph where loops are not allowed. Every function of an empty digraph is a homomorphism, so these are quasi-projective digraphs. A homomorphism of a complete digraph (to any loopless digraph) is injective. Thus $j$ is an isomorphism. Thus $\phi=j^{-1}f$ satisfies the required condition.

 If the digraph, where loops are not allowed, is quasi-projective, then the proof of proving that it is an empty or a complete digraph is similar as in the proof of Lemma \ref{g1}. We may assume that $(v_1, v_2) \notin A$ and $(v_3 v_4)\in A$ and then we map the digraph to the a complete digraph by the maps $f$ and $j$ in the same way as in the proof.
\end{proof}

Now we characterise graphs and digraphs with loops.

\begin{lemma}\label{gl1}
In a non-empty quasi-projective graph, where loops are allowed, every vertex has a loop.
\end{lemma}

\begin{proof} Let a non-empty graph ${\cal G}=(V,E)$ be quasi-projective, where $V=\{v_1, v_2,\ldots, v_n\}$ and let $(v_1,v_2)$ be an edge. Suppose to the contrary that there is a vertex $v$ that it does not have a loop. Let $T=(\{k_1,k_2\},E_T)$, where $E_T=\{(k_1,k_1),(k_1,k_2),(k_2,k_2)\}$. We define $f$ and $j$ in the following way:

$$f(x)=\begin{cases}   k_1 & if \hbox{ } x=v_i \end{cases},$$

$$j(x)=\begin{cases}  k_{1} & if \hbox{ } x=v \\ k_2 & otherwise \end{cases}.$$

The map $f$ is homomorphism and the map $j$ is epimorphism obviously. Now we prove, there is no homomorphism $\phi$ such that $j \circ \phi=f$.  We have $f(v_1)=f(v_2)=k_1$, the preimage of $k_1$ is $\{v\}$ by the map $j$. So, if there is $\phi$ such that $f=j \circ \phi$, then $\phi(v_1)=v$ and $\phi(v_2)=v$. This contradicts the fact that $(v,v)\notin E$.

\end{proof}

\begin{lemma}\label{gl2}
In an e-empty quasi-projective graph, where loops are allowed, if a vertex has a loop, then every vertex has a loop.
\end{lemma}

\begin{proof} Let an empty graph ${\cal G}=(V,E)$ be quasi-projective, where $V=\{v_1, v_2,\ldots, v_n\}$ and let $(v_1,v_1)$ be a loop. Suppose to the contrary that there is a vertex $v$ that it does not have a loop. Let $T=(\{k_1,k_2\},E_T)$, where $E_T=\{(k_1,k_1),(k_1,k_2),(k_2,k_2)\}$. We define the same $f$ and $j$ functions as in Lemma \ref{gl1}.

 Now we prove, there is no homomorphism $\phi$ such that $j \circ \phi=f$.  We have $f(v_1)=k_1$, the preimage of $k_1$ is $\{v\}$ by the map $j$. So, if there is $\phi$ such that $f=j \circ \phi$, then $\phi(v_1)=v$. This contradicts the fact that $(v,v)\notin E$.

\end{proof}

\begin{lemma}\label{gl3}
A quasi-projective graph, where loops are allowed, is a complete graph or an e-empty graph where if there is a loop, every vertex has a loop.
\end{lemma}

\begin{proof}
It follows from Lemmata \ref{g1}, \ref{gl1} and \ref{gl2}.

\end{proof}

\begin{lemma}\label{gl4}
Let ${\cal G}=(V,E)$ be a graph where loops are allowed. If it is a complete graph or an e-empty graph where if there is a loop, every vertex has a loop, then it is quasi-projective.
\end{lemma}

\begin{proof} Let ${\cal G}=(V,E)$ be an e-empty graph where if there is a loop, every vertex has a loop. Every function of an e-empty graph without loops or when every vertex has a loop is a homomorphism, so these are quasi-projective graphs. Now let ${\cal G}=(V,E)$ be a complete graph. Let $f: {\cal G}\rightarrow {\cal H}$ and $j: {\cal G}\rightarrow {\cal H}$ be homomorphism and epimorphism, respectively, where ${\cal H}=(V_H,E_H)$ is a graph. We define $\phi$ in the following way: we map vertex $v\in V$ to a vertex from $j^{-1}(f(v))$. Because $\cal G$ is a complete graph where every vertex has a loop, then $f({\cal G})$ is a complete graph where every vertex has a loop. Also, $j^{-1}(f({\cal G}))$ is a complete graph where every vertex has a loop, thus $\phi$ is a homomorphism.

\end{proof}

\begin{theorem}\label{gl5}
A digraph, where loops are allowed, is quasi-projective if and only it is a complete graph or an e-empty graph where if there is a loop, every vertex has a loop.
\end{theorem}

\begin{proof}
Let  ${\cal D}=(V,A)$ be an e-empty digraph where if there is a loop, every vertex has a loop. Every function of an empty digraph without loops or when every vertex has a loop is a homomorphism, so these are quasi-projective graphs. Now let ${\cal D}=(V,A)$ be a complete digraph where every vertex has a loop. Let $f: {\cal D}\rightarrow {\cal H}$ and $j: {\cal D}\rightarrow {\cal H}$ be homomorphism and epimorphism, respectively, where ${\cal H}=(V_H,E_H)$ is a digraph. We define $\phi$ in the following way: we map vertex $v\in V$ to a vertex from $j^{-1}(f(v))$. Because $\cal D$ is a complete digraph where every vertex has a loop, then $f({\cal D})$ is a complete digraph where every vertex has a loop. Also, $j^{-1}(f({\cal D}))$ is a complete digraph where every vertex has a loop, thus $\phi$ is a homomorphism.

 If the digraph is quasi-projective, then the proof of proving that it is a complete graph or an e-empty graph where if there is a loop, every vertex has a loop, is similar as in the proof of Lemma \ref{gl3}.

\end{proof}

\section{Quasi-projective  hypergraphs} \label{qphg}

In this section we characterise quasi-projective finite hypergraphs. 

A {\it hypergraph} ${\cal H}$ is a pair ${\cal H}=(V,E)$ where $X$ is a set of elements called {\it nodes} or {\it vertices}, and $E$ is a set of non-empty subsets of $V$ called {\it hyperedges} or {\it edges}. Therefore, $E$ is a subset of ${\mathcal {P}}(V)\setminus \{\emptyset \}$, where ${\mathcal {P}}(V)$ is the power set of $V$. The size of the vertex set is called the {\it order} of the hypergraph, and the size of edges set is the {\it size} of the hypergraph, we denote this by $n({\cal H})$ and $m({\cal H})$, respectively.

We say that a hypergraph is {\it $k$-uniform} if each hyperedge has size $k$. In particular, the {\it complete $k$-uniform hypergraph} on  vertices has all $k$-subsets of $n({\cal H})$ as edges. We denote this by $K_n^k$.

A mapping $f \colon V_1 \rightarrow V_2$ is a {\it homomorphism} from a hypergraph ${\cal H}_1=(V_1,E_1)$ to a hypergraph ${\cal H}_2=(V_2,E_2)$ if $e \in E_1$, then the set of images of all vertices from the set $e$ by the map $f$ is in $E_2$.

\begin{theorem} \label{hg1}
If a hypergraph is quasi-projective, then if there is an edge with $n$ vertices, then the set of every at most $n$ vertices is an edge.
\end{theorem}

\begin{proof}
Let ${\cal H}=(V,E)$ be a quasi-projective hypergraph, where $|V|=n$. Suppose that there is an edge $e_1=\{v_1,v_2, \ldots, v_k\}$, and assume $\{u_1, \ldots u_l \}$ is not an edge, where $l \leq k$. Let ${\cal T}=(V_T,E_T)$ be a hypergraph, where $V_T=\{v^T_1, v_2^T, \ldots, v^T_{l+1}\}$ and $E_T={\cal P}(V_T)\setminus \{\emptyset\}$. We define $f$ and $j$ in the following way:

$$f(x)=\begin{cases}  v_i^T & if \hbox{ }  x=v_i,\hbox{ } i=1, 2,\ldots,l \\ v_l^T &   if \hbox{ }  x=v_i,\hbox{ } i=l+1, l+2,\ldots,k \\  v^T_{l+1} & otherwise \end{cases}$$

$$j(x)=\begin{cases} v_i^T & if \hbox{ }  x=u_i, \hbox{ } i=1,2,\ldots,l \\
v^T_{l+1} &  otherwise \end{cases}$$

The map $f$ is homomorphism and the map $j$ is epimorphism obviously. We have $f(\{v_i\}_{i\in\{1,2,\ldots,k\}})=\{v^T_i\}_{i\in\{1,2,\ldots,l\}}$ and the set of $j^{-1}(\{v_i^T\}_{i\in \{1,2,\ldots,l\}})$ is not an edge. Thus does not exist $\phi$ such that $f=j \circ \phi$.

\end{proof}

\begin{theorem} \label{hg2}
A hypergraph is quasi-projective if and only if there is an edge with $n$ vertices, then the set of every at most $n$ vertices is an edge.
\end{theorem}

\begin{proof} If there is a quasi-projective hypergraph, then from Lemma \ref{hg1} we know that if there is an edge with $n$ vertices, then the set of every at most $n$ vertices is an edge.

Let ${\cal H}=(V,E)$ be a hypergraph where if there is an edge with $n$ vertices, then the set of every at most $n$ vertices is an edge. Every function of an empty graph is a homomorphism, so these are quasi-projective hypergraphs. Because every function from a hypergraph where if there is an edge with $n$ vertices, then the set of every at most $n$ vertices is an edge, into itself is a homomorphism, it is quasi-projective.

\end{proof}

\section{Quasi-projective point-line geometries}\label{qpplg}

In this section we characterise quasi-projective finite point-line geometries. 

A {\it point-line geometry} is an ordered pair $(\cal X,L)$, where $\cal X$ is a non-empty set of elements called
{\it points} and $\cal L \subseteq P(X)$ is a collection of subsets called {\it lines} such that  every line contains at least two points and every pair of distinct points is contained in at most one line. We only consider finite point-line geometries.

A {\it subgeometry or substructure} $(\cal Y,L_Y)$ of the point-line geometry $(\cal X,L)$ is a point-line geometry, where
$\emptyset \not= \cal Y \subseteq X$ and ${\cal L_Y}=\{l \cap{\cal  Y} \mid l \in {\cal L}
\wedge  |l\cap {\cal Y}| \geq 2 \}$. If $\emptyset\not= {\cal Y}\subseteq {\cal X}$ then the point-line geometry $(\cal Y,L_Y)$ induced on $\cal Y$,  where ${\cal L_Y}=\{l \cap{\cal  Y} \mid (\forall l) (l \in {\cal L}
\wedge  |l\cap {\cal Y}| \geq 2 )\}$, is an induced subgeometry of $({\cal X,\cal L})$.

A line which contains more than two points is called a {\it regular line}. A line which contains exactly two points is called a {\it singular line}. Regular lines will be denoted by lower case letters $a,b,c, \dots$ and singular lines will mostly be denoted as $AB$, where $A$ and $B$ are the points contained in it. An isolated point is a point which belongs to no line of the geometry. The points $A$ and $B$ are collinear if there exists a line $l\in {\cal L}$ such that $A,B\in l$. In this case we write $A \sim B$. 

A mapping $f : \cal X \rightarrow Y$ is a {\it homomorphism} from a point-line geometry $(\cal X, L)$ to a point-line geometry $\cal (Y,K)$ if for every $l \in \cal L$, either $|f(l)|=1$ or there is a line $k \in \cal K$ such that $f(l) \subseteq k$. A homomorphism $f$ from an arbitrary induced subgeometry of a point-line geometry $(\cal X,L)$ into $(\cal X,L)$ will be referred to as a {\it  local homomorphism} of $(\cal X,L)$. 

The point-line geometries in our paper are not first-order structures. In order to relate them to the results mentioned in the introduction, see \cite{jungabel1} where authors define point-line geometries as first-order structures and show that the homomorphisms and substructures of this structure are canonically the same as the homomorphisms and subgeometries according to our definitions. This way we can consider the class of point-line geometries to be a subclass of 3-uniform hypergraphs.

\begin{lemma}\label{plg1}
Let $({\cal X},{\cal L})$ be a quasi-projective point-line geometry. Then it is a geometry:

\begin{enumerate}

\item  without lines, or

\item where every two points are on singular lines, or

\item where every point is on a regular line.

\end{enumerate} 
\end{lemma}

\begin{proof} 
Let $({\cal X},{\cal L})$ be a quasi-projective point-line geometry, where $|{\cal X}|=n$. Suppose that there is a singular line and there are two points which are not on a singular line. Then the proof is the same as in Lemma \ref{g1} where there is a graph. So, now suppose that there is a line $r_1$ with $r_1=\{A_1,A_2, \ldots, A_k\}$, $k\geq 3$ and let $B_1, B_2,\ldots, B_{k}$ be points such that they are not on a line. Let $({\cal T},{\cal L}_T)$ be a point-line geometry, where ${\cal T}=\{A^T_1, A_2^T, \ldots, A^T_k\}$ and ${\cal L}_T=\{r_T\}$ and $r_T=\{A^T_1, A_2^T, \ldots, A^T_k\}$. We define $f$ and $j$ in the following way:

$$f(X)=\begin{cases}  A_i^T & if \hbox{ }  X=A_i , \hbox{ } i\in\{1,2,\ldots,k\}  \\  A^T_k & otherwise \end{cases}$$

$$j(X)=\begin{cases}  A^T_i & if \hbox{ }  X=B_i, \hbox{ } i\in\{1,2,\ldots, k\} \\  A^T_2 & if \hbox{ }  X \hbox{ \it form a line with $B_2, B_3,\ldots,B_{k}$} \\ A_1^T & otherwise. \end{cases}$$

The map $f$ is homomorphism and the map $j$ is epimorphism obviously. We have $f(A_i)=A^T_i$, for $i\in \{1,2,\ldots,k\}$. The set  $j^{-1}(A_i^T)$, for $i\in \{1,2,\ldots,k\}$ is $\cal X$ and there is not a line $r$ such that $j(r)=r_T$. If we observe points $B_1,B_2,\ldots,B_{k}$, then they are not on a line. If we consider points $B_1,X,B_3,\ldots,B_{k}$, where points $X,B_2,B_3,\ldots,B_{k}$ are on a line, then if we suppose that points $B_1,X,B_3,\ldots,B_{k}$ are on a line, then points $B_1,B_2,B_3,\ldots, B_{k}$ are on a line which is a contradiction. If we look at points $Y,B_2,B_3,\ldots,B_{k}$, then these points cannot form a line unless from definition of the function $j$, then $\phi(Y)=A^T_2$. Also, if we observe points $Y,X,B_3,\ldots,B_{k}$, where points $X,B_2,B_3,\ldots,B_{k}$ form a line, then these points cannot form a line unless points $Y,B_2,B_3,\ldots,B_{k}$ form a line and from definition of the function $j$, then $\phi(Y)=A^T_2$.

So, we have got that if there is a line with $k$ points, then every $k$ points are on a line. If $k=2$ or $k=n$, then every two points are on singular lines or every point is on a regular line, respectively. The case when $2<k<n$ leads to a contradiction, because then there are two points which are on more than one line.

\end{proof}

\begin{theorem} \label{plg2}
A point-line geometry is quasi-projective if and only if it is a geometry:

\begin{enumerate}

\item  without lines, or

\item where every two points are on singular lines, or

\item where every point is on a regular line.

\end{enumerate} 
\end{theorem}

\begin{proof} If there is a quasi-projective geometry, then from Lemma \ref{plg1} we know that it is a geometry without lines or a geometry where every two points are on singular lines or a geometry where every point is on a regular line.

Let $({\cal X},{\cal L})$ be a point-line geometry. Every function of a geometry without line is a homomorphism, so these are quasi-projective geometries. Because every function from a geometry, where every two points are on singular lines or where every point is on a regular line, into itself is a homomorphism, it is quasi-projective.

\end{proof}

\noindent
{\bf Acknowledgement.} The author would like to express gratitude to Csaba Szab\'o  and G\'abor Somlai  for their careful reading, comments and remarks. They highly improved the presentation of this article.

\bibliographystyle{plain}

\vskip5mm
\noindent
{\sc \' Eva Jung\'abel, E\"otv\"os Lor\'and University, P\'azm\'any P\'eter s\'et\'any 1/C, 1117 Budapest, Hungary}

\noindent
{\it E-mail address:} {\tt evajungabel@student.elte.hu}

 \end{document}